\documentclass[10pt]{article}
\usepackage{latexsym}
\usepackage{amsbsy}
\usepackage{amsfonts,mathrsfs,amssymb}
\usepackage{amsmath,amsthm}
\usepackage[latin1]{inputenc}
\usepackage{graphics,graphicx}
\usepackage{a4wide}

 \newtheorem{thm}{Theorem}[section]
 \newtheorem{defin}[thm]{Definition}
 \newtheorem{lem}[thm]{Lemma}
 \newtheorem{prop}[thm]{Proposition}
 \newtheorem{cor}[thm]{Corollary}
 \newtheorem{rem}[thm]{Remark}
 \newtheorem{ex}[thm]{Example}

 \newcommand{\bthm}{\begin{thm}}
 \newcommand{\ethm}{\end{thm}}
 \newcommand{\bd}{\begin{defin}}
 \newcommand{\ed}{\end{defin}}
 \newcommand{\blem}{\begin{lem}}
 \newcommand{\elem}{\end{lem}}
 \newcommand{\bcor}{\begin{cor}}
 \newcommand{\ecor}{\end{cor}}
 \newcommand{\bprop}{\begin{prop}}
 \newcommand{\eprop}{\end{prop}}
 \newcommand{\brem}{\begin{rem} \rm}
 \newcommand{\erem}{\end{rem}}
 \newcommand{\bex}{\begin{ex} \rm}
 \newcommand{\eex}{\end{ex}}

 \newcommand{\pr}{\noindent{\bf Proof. }}
 \newcommand{\ep}{\nolinebreak{\hspace*{\fill}$\Box$ \vspace*{0.25cm}}}

 \newcommand{\beq}{\begin{equation} }
 \newcommand{\eeq}{\end{equation} }

 \newcommand{\bea}{\begin{eqnarray}}
 \newcommand{\eea}{\end{eqnarray}}

 \newcommand{\beas}{\begin{eqnarray*}}
 \newcommand{\eeas}{\end{eqnarray*}}

 \newcommand{\beqs}{\begin{equation*}}
 \newcommand{\eeqs}{\end{equation*}}

 \newcommand{\bite}{\begin{itemize}}
 \newcommand{\eite}{\end{itemize}}

 \newcommand{\ben}{\begin{enumerate}}
 \newcommand{\een}{\end{enumerate}}

 \newcommand{\ba}{\begin{array}}
 \newcommand{\ea}{\end{array}}

 \newcommand{\ds}{\displaystyle}

 \newcommand{\R}{\mathbb R}
 \newcommand{\N}{\mathbb N}

 \newcommand{\cA}{\ensuremath{{\cal A}}}
 \newcommand{\cC}{\ensuremath{{\cal C}}}

 \newcommand{\cL}{\ensuremath{{\cal L}}}

 \newcommand{\cS}{\ensuremath{{\cal S}}}

 \newcommand{\pd}{\partial}

 \newcommand{\eps}{\varepsilon}
 \newcommand{\vphi}{\varphi}

 \newcommand{\supp}{\ensuremath{\mathrm{supp\,}}}

\begin{document}

 \title{
 Variational Problems with Fractional Derivatives: Euler-Lagrange Equations
 }

 \author{
 Teodor M. Atanackovi\'c
 \footnote{Faculty of Technical Sciences, Institute of Mechanics, University of Novi Sad, Trg Dositeja Obradovi\' ca 6, 21000 Novi Sad, Serbia.
 Electronic mail: atanackovic@uns.ac.rs.
 Work supported by
 Projects 144016 and 144019 of the Serbian Ministry of Science
 and
 START-project Y-237 of the Austrian Science Fund.}\\
 Sanja Konjik
 \footnote{Faculty of Agriculture, Department of Agricultural Engineering, University of Novi Sad, Trg Dositeja Obradovi\' ca 8, 21000 Novi Sad, Serbia.
 Electronic mail: sanja\_konjik@uns.ac.rs. }\\
 Stevan Pilipovi\'c
 \footnote{Faculty of Sciences, Department of Mathematics, University of Novi Sad, Trg Dositeja Obradovi\' ca 4, 21000 Novi Sad, Serbia.
 Electronic mail: pilipovic@dmi.uns.ac.rs}\\
 }

 \date{}
 \maketitle

 \begin{abstract}
 We generalize the fractional variational problem by allowing the possibility
 that the lower bound in the fractional derivative does not coincide
 with the lower bound of the integral that is minimized. Also, for
 the standard case when these two bounds coincide, we derive a new
 form of Euler-Lagrange equations. We use approximations for
 fractional derivatives in the Lagrangian and obtain the Euler-Lagrange
 equations which approximate the initial Euler-Lagrange equations
 in a weak sense.

 \vskip5pt


 \vskip5pt


 \vskip5pt

 \end{abstract}

 \section{Introduction}

 Fractional calculus with derivatives and integrals of any real or
 complex order has its origin in the work of L. Euler, and even
 earlier in the work of G. Leibnitz.
 Shortly after being introduced, the new theory turned out to be
 very attractive to many famous mathematicians and scientists
 (e.g. P. S. Laplace, B. Riemann,
 J. Liouville, N. H. Abel, J. B. J. Fourier, et al.) due to the
 numerous possibilities for its applications. Beside in mathematics,
 fractional derivatives and integrals appear in
 physics, mechanics, engineering, elasticity, dynamics, control
 theory, electronics, modelling, probability, finance, economics,
 biology, chemistry, etc. The fractional
 calculus is nowadays covered with several extensive reference
 books \cite{GorenfloMainardi, KilSriTru, MillRoss, OldSpa, SamkoKM}
 and a large number of relevant papers.

 Fractional calculus of variations unifies calculus of variations
 (cf.\ classical books \cite{GelfandF, GiaquintaH, Jost, Sagan, vanBrunt})
 and fractional calculus, by inserting fractional derivatives into
 variational integrals. This of course occurs naturally in many
 problems of physics or mechanics, in order to provide more
 accurate models of physical phenomena.
 Research within this topic goes in different directions.
 G. Jumarie \cite{Jumarie93, Jumarie95} is one of the first who
 has used fractional variational calculus in the analysis of
 fractional Brownian motion. Especially, we refer to
 G. Jumarie's paper \cite{Jumarie06} for the new approach to
 fractional stochastic mechanics and stochastic optimal control.
 We cite also F. Riewe \cite{Riewe96, Riewe97}, who investigated
 nonconservative Lagrangian and Hamiltonian mechanics and
 for those cases formulated a version of the Euler-Lagrange
 equations. Further study of the fractional Euler-Lagrange
 equations can be found in the work
 of O. P. Agrawal \cite{Agrawal02, Agrawal06, Agrawal07}, who
 considered different types of variational problems,
 involving Riemann-Liouville, Caputo and Riesz fractional
 derivatives, respectively. He derived corresponding Euler-Lagrange
 equations and discussed possibilities for prescribing
 boundary conditions in each case. The work of mentioned
 authors influenced many recent
 papers. For instance, D. Baleanu \cite{Bal04, Bal05} applied the
 fractional Euler-Lagrange equations to examine fractional
 Lagrangian and Hamiltonian systems linear in velocities.
 Other applications of fractional variational principles are
 presented in \cite{AOP07, AtanackovicStankovic07, Bal06,
 DreisigmeyerYoung04, DreisigmeyerYoung03, GaiesElAkrmi,
 Klimek01, Lazopoulos06, Odibat07, Odibat06, Baleanu&co07}.
 We cite here also the work of
 G. S. F. Frederico and D. F. M. Torres \cite{FredericoTorres},
 who introduced a new concept of fractional conserved quantities
 on the basis of a variational principle, and proved a
 version of the fractional N\" other theorem.

 There are several aims of this paper. First, we discuss the
 Euler-Lagrange equations of \cite{Agrawal02} and \cite{Agrawal06}
 and show that the transversality condition proposed in \cite{Agrawal06}
 should be used with care, since it may lead to erroneous conclusions.
 Second, we consider a fractional variational problem, defined
 by a functional whose lower bound does not coincide with the lower
 bound in the left Riemann-Liouville fractional derivative that
 appears in the Lagrangian. This leads to the natural generalization
 of the fractional variational problems considered so far.
 Third, in Section \ref{sec:approximations} we approximate
 fractional derivative in the Lagrangian $L(t,u(t), {}_aD_t^\alpha u)$
 with a finite number of terms containing derivatives
 of integer order, which reduces the variational problem to the
 one depending only on the classical derivatives of the
 function $u$. For that purpose we consider approximations in
 a weak sense, using analytic functions as a test
 function space, and show that a sequence of approximated
 Euler-Lagrange equations converges to the fractional
 Euler-Lagrange equation.

 \section{Notation}

 Let $u\in L^1([a,b])$ and $0\leq\alpha,\beta <1$. Then the left
 Riemann-Liouville fractional integral of order $\alpha$,
 ${}_aI_t^\alpha u$ is defined as
 $$
 {}_aI_t^\alpha u=\frac{1}{\Gamma (\alpha)}
 \int_a^t (t-\theta)^{\alpha -1} u(\theta)\, d\theta,
 \quad t\in [a,b].
 $$
 The right Riemann-Liouville fractional integral of order
 $\beta$, ${}_tI_b^\beta u$ is defined as
 $$
 {}_tI_b^\beta u=\frac{1}{\Gamma (\beta)}
 \int_t^b (\theta-t)^{\beta -1} u(\theta)\, d\theta,
 \quad t\in [a,b].
 $$
 If $u$ is an absolutely continuous function in $[a,b]$, i.e.
 $u\in AC([a,b])$, and $0\leq\alpha <1$, then the left
 Riemann-Liouville fractional derivative of order $\alpha$,
 ${}_aD_t^\alpha u$ is given by
 $$
 {}_aD_t^\alpha u=\frac{d}{dt} {}_aI_t^{1-\alpha} u=
 \frac{1}{\Gamma (1-\alpha)}\frac{d}{dt}
 \int_a^t (t-\theta)^{-\alpha} u(\theta)\, d\theta,
 \quad t\in [a,b],
 $$
 and for $0\leq\beta <1$,
 the right Riemann-Liouville fractional derivative of order
 $\beta$, ${}_tD_b^\beta u$ is given by
 $$
 {}_tD_b^\beta u=\Big(-\frac{d}{dt}\Big) {}_tI_b^{1-\beta} u=
 \frac{1}{\Gamma (1-\beta)}\Big(-\frac{d}{dt}\Big)
 \int_t^b (\theta-t)^{-\beta} u(\theta)\, d\theta,
 \quad t\in [a,b].
 $$
 We have that ${}_aD_t^\alpha {}_aI_t^\alpha =I$, where $I$ is the identity map.
 The formula for fractional integration by parts reads
 (see \cite[p.\ 46]{SamkoKM}):
 \beq \label{eq:frac int by parts}
 \int_a^b f(t) {}_aD_t^\alpha g\, dt = \int_a^b g(t)
 {}_tD_b^\alpha f\, dt, \quad f,g\in AC([a,b]).
 \eeq

 In the distributional setting, the Riemann-Liouville fractional
 derivatives can be defined via convolutions in the space of
 tempered distributions supported by $[0,+\infty)$.
 Let
 $(f_\alpha)_{\alpha\in\R}\in\cS'_{+}=\{f\in\cS'(\R):\supp f\subset
 [0,+\infty)\}$ be a family of distributions defined as
 $$
 f_\alpha(t)=\left\{
 \ba{ll}
 \ds H(t)\frac{t^{\alpha-1}}{\Gamma(\alpha)}, & \alpha>0\\
 f_{\alpha+N}^{(N)}(t), & N\in\N: N+\alpha>0 \wedge N+\alpha-1<0
 \ea
 \right.
 $$
 where $H$ is the Heaviside function. The operator
 $f_\alpha\ast$ is the Riemann-Liouville operator of
 differentiation, resp. integration of order $\alpha$
 for $\alpha<0$, resp. $\alpha>0$. In this setting ${}_aI_t^\alpha$ and
 ${}_aD_t^\alpha$ are inverses in both directions, due to the group
 property in $\cS'$, that is $f_\alpha\ast f_\beta=f_{\alpha+\beta}$,
 for all $\alpha,\beta\in\R$.

 Beside the Riemann-Liouville approach there exist
 several other possibilities for introducing derivatives of
 fractional order. We will make use of the Caputo
 fractional derivatives: if $0\leq\alpha,\beta <1$ and
 $u\in AC([a,b])$, then the left Caputo fractional
 derivative of order $\alpha$, ${}_a^cD_t^\alpha u$
 is defined as
 $$
 {}_a^cD_t^\alpha u= \frac{1}{\Gamma (1-\alpha)}
 \int_a^t (t-\theta)^{-\alpha} \dot{u}(\theta)\, d\theta,
 \quad t\in [a,b],
 $$
 where $\dot{(\cdot)}$ denotes the total derivative
 $\ds \frac{d}{dt}$,
 and the right Caputo fractional
 derivative of order $\beta$, ${}_t^cD_b^\beta u$
 is defined as
 $$
 {}_t^cD_b^\beta u= - \frac{1}{\Gamma (1-\beta)}
 \int_t^b (\theta-t)^{-\beta} \dot{u}(\theta)\, d\theta,
 \quad t\in [a,b].
 $$
 The Riemann-Liouville and Caputo fractional derivatives are
 related by the following formula:
 \beq \label{eq:veza RL i Caputo-levi}
 {}_aD_t^\alpha u={}_a^cD_t^\alpha u + \frac{1}{\Gamma (1-\alpha)}
 \frac{u(a)}{(t-a)^\alpha},
 \quad t\in [a,b],
 \eeq
 and similarly
 \beq \label{eq:veza RL i Caputo-desni}
 {}_tD_b^\beta u={}_t^cD_b^\beta u + \frac{1}{\Gamma (1-\beta)}
 \frac{u(b)}{(b-t)^\alpha},
 \quad t\in [a,b].
 \eeq

 If $u(a)=0$ (resp.\ $u(b)=0$) then the left (resp.\ right)
 Riemann-Liouville and Caputo fractional derivatives coincide.
 Also, for $u(a)=0$ the left Riemann-Liouville fractional
 derivative commutes with the first derivative with
 respect to $t$, i.e. $\frac{d}{dt}{}_aD_t^\alpha u=
 {}_aD_t^\alpha \frac{d}{dt} u$ (and the same holds for the
 right Riemann-Liouville fractional derivative if $u(b)=0$).

 In this paper we will consider the Lagrangian $L$ as
 a function of $t,u$ and ${}_aD_t^\alpha u$, i.e.\ $L=L(t,u(t),
 {}_aD_t^\alpha u)$. The partial derivatives of $L$ will be
 denoted by $\ds \frac{\pd L}{\pd t}$, $\ds \frac{\pd L}{\pd u}$ and
 $\ds \frac{\pd L}{\pd {}_aD_t^\alpha u}$, or by $\pd_1 L$, $\pd_2 L$
 and $\pd_3 L$ respectively. The first
 (or Lagrangian) variation will be denoted by $\delta$, as usual.

 \section{Euler-Lagrange equations}

 Let $(A,B)$ be a subinterval of $(a,b)$. Consider a functional
 \beq \label{eq:fractional vp}
 \cL[u]=\int_A^B L(t,u(t),{}_aD_t^\alpha u)\, dt,
 \quad 0\leq \alpha< 1,
 \eeq
 where $u$ is an absolutely continuous function in $[a,b]$ and
 $L$ is a function in $(a,b)\times \R\times \R$ such that
 \beq \label{eq:conditions on L}
 \left.\ba{c}
 L\in\cC^1((a,b)\times \R\times \R)\\
 \mbox{ and }\\
 t\mapsto\pd_3 L(t,u(t),{}_aD_t^\alpha u)
 \in AC([a,b]), \mbox{ for every } u\in AC([a,b])
 \ea\right\}
 \eeq
 A fractional variational problem consists of finding
 extremal values (minima or maxima) of the functional
 (\ref{eq:fractional vp}) among all admissible functions.
 The function $L$ is called the Lagrangian. Note that in
 (\ref{eq:fractional vp}) the constants $a$ and $A$ are
 assumed to be different and in general and $a\leq A$.
 Their physical meaning is
 also different. While the interval $(A,B)$ defines
 the Hamilton action, the value $a$ defines memory of the system.
 In the special case, that was treated previously (cf.\
 \cite{Agrawal02, Agrawal06, Bal04, Bal05, FredericoTorres}),
 it was assumed that $a=A$.

 As mentioned above, $\cL$ is to be minimized (or maximized) over
 the set of admissible functions. Hence, we have to specify where
 we look for a minimum of (\ref{eq:fractional vp}): the admissible
 set will consists of all absolutely continuous functions
 $u$ in $[a,b]$, which pass
 through a fixed point at $a$, i.e. $u(a)=a_0$, for a fixed
 $a_0\in\R$.

 \brem
 (i)\ We consider a fractional variational problem which
 involves only left Riemann-Liouville fractional derivatives. The
 problem can be easily generalized for Lagrangians which will
 depend also on the right Riemann-Liouville fractional
 derivatives.\\
 (ii)\ We assume that $0\leq\alpha<1$. Our assumption can
 be extended to the case $\alpha\geq 1$ without difficulties.\\
 (iii)\ Traditionally, the minimizers of a variational problem are
 sought. Analogously, one can consider the problem of
 finding the maximal values of a variational problem.
 \erem

 In this section we discuss the results on the Euler-Lagrange
 equations of the fractional variational problem (\ref{eq:fractional
 vp}).

 The case $A=a$ was treated in \cite{Agrawal02} by Agrawal.
 It was proved there that if one wants to minimize
 (\ref{eq:fractional vp}) among all functions $u$ which have
 continuous left $\alpha$-th Riemann-Liouville fractional
 derivative and which satisfy the Dirichlet boundary conditions
 $u(a)=a_0$ and $u(b)=b_0$, for some real constant values $a_0$ and
 $b_0$, then a minimizer should be sought among all solutions
 of the Euler-Lagrange equation
 \beq \label{eq:EL with RL}
 \frac{\pd L}{\pd u} + {}_tD_b^\alpha \Big(\frac{\pd L}{\pd _aD_t^\alpha
 u}\Big) =0.
 \eeq
 This result was modified in \cite{Agrawal06}, where again $a=A$
 was used, and the boundary condition
 was specified at $t=a$ only, which allowed the natural boundary
 conditions to be developed. The corresponding Euler-Lagrange
 equation was obtained as
 \beq \label{eq:EL with Caputo}
 \frac{\pd L}{\pd u} + {}_t^cD_b^\alpha \Big(\frac{\pd L}{\pd _aD_t^\alpha
 u}\Big) =0,
 \eeq
 with the transversality (natural) condition
 \beq \label{eq:transversality cond}
 \frac{\pd L}{\pd _aD_t^\alpha u} {}_aI_t^{1-\alpha}\delta u =0
 \quad \mbox{ at } t=b.
 \eeq
 It is clear that (\ref{eq:EL with Caputo}) and
 (\ref{eq:transversality cond}) imply (\ref{eq:EL with RL}). But
 the converse does not hold in general. This depends on
 assumptions on $L$ and the set of admissible functions. Our first
 example which is to follow will show that assuming
 (\ref{eq:conditions on L}), condition (\ref{eq:EL with RL}) does
 not imply (\ref{eq:EL with Caputo}) and
 (\ref{eq:transversality cond}).

 \bex
 Consider a fractional variational problem of the form
 \beq \label{eq:frac vp 01}
 \cL[u]=\int_0^1 L(t,u(t),{}_0D_t^\alpha u)\, dt \to\min
 \eeq
 with $u\in AC([0,1])$ and $L$ to be specified.\\
 (i)\ Let the Lagrangian $L$ satisfies (\ref{eq:conditions on L})
 and has the form
 $$
 L(t,u,{}_aD_t^\alpha u)= F(t,u)+f(t){}_aD_t^\alpha u.
 $$
 Then $\pd L/\pd {}_aD_t^\alpha u=f(t)$, $t\in (0,1)$, and we have
 (see \cite[(13)]{Agrawal06}) that
 \beas
 0 &=& \int_0^1 \Big(\frac{\pd F}{\pd u}+{}_t^cD_1^\alpha f\Big)\delta
 u\, dt + f(t){}_0D_t^{\alpha-1} \delta u {\Big\vert}_{_{0}}^{^{1}}
 \\
 &=& \int_0^1 \Big(\frac{\pd F}{\pd u}
 + \frac{1}{\Gamma
 (1-\alpha)}\int_t^1\frac{f'(\theta)}{(\theta-t)^\alpha}\, d\theta
 + \frac{f(1)}{\Gamma(1-\alpha)}\frac{1}{(1-t)^{\alpha}}
 \Big)\delta u\, dt.
 \eeas
 Now, take for example $F=\frac{u^2}{2\Gamma(1-\alpha)}
 \frac{1}{(1-t)^{\alpha}}$ and $f(t)=-1, t\in (0,1)$. Then the
 Euler-Lagrange equation (\ref{eq:EL with RL}) gives
 $$
 \frac{u}{\Gamma(1-\alpha)}\frac{1}{(1-t)^{\alpha}}=
 \frac{1}{\Gamma(1-\alpha)}\frac{1}{(1-t)^{\alpha}},
 $$
 thus $u\equiv 1$ in $[0,1]$.
 Hence, if one wants to formulate (\ref{eq:frac vp 01}) so that
 (\ref{eq:EL with RL}) holds if and only if (\ref{eq:EL with Caputo})
 and (\ref{eq:transversality cond}) hold, then
 some additional assumptions on $F$ have to be supposed.
 For instance, the condition $f(1)=0$ provides the desired equivalence
 of (\ref{eq:EL with RL}) with (\ref{eq:EL with Caputo},
 \ref{eq:transversality cond}).
 \\
 (ii)\ The Lagrangian
 $$
 L(t,u,{}_aD_t^\alpha u)=({}_aD_t^\alpha u-u)^2
 $$
 was investigated in details in \cite{Agrawal06}.
 Solutions of the corresponding fractional variational problem
 (\ref{eq:frac vp 01}) can be found directly. It is clear that the
 functional $\cL$ achieves its minimum
 (which is zero) when ${}_aD_t^\alpha u-u=0$. Hence, the problem reduces
 to solving the equation ${}_aD_t^\alpha u=u$.
 It was shown in \cite[p.\ 222]{KilSriTru} that this equation has no
 solution which is bounded at $0$. In other
 words, one can not solve the fractional variational problem with the
 Lagrangian given above among the functions with the prescribed,
 finite boundary condition at $0$.
 If instead of the left Riemann-Liouville fractional derivative
 one considers the Lagrangian $L$ as a function of the left Caputo
 fractional derivative, then the equation ${}_a^cD_t^\alpha u=u$
 has a solution bounded at zero, which is also a solution of the
 corresponding fractional variational problem (\ref{eq:frac vp 01}).
 In a recent paper \cite{Agrawal07}, the author
 tried to overcome the problem of non-solvability of the equation
 ${}_aD_t^\alpha u=u$ by using the
 symmetrized Caputo fractional derivative, called the Riesz Caputo
 fractional derivative and defined as
 ${}_a^{rc}D_t^\alpha u:=\frac{1}{2} ({}_a^cD_t^\alpha u-
 {}_t^cD_b^\alpha u)$. Similar kind of fractional derivatives were
 used earlier in \cite{Klimek01, Lazopoulos06}.
 \eex

 These examples suggest that the Euler-Lagrange equation for the
 fractional variational problem (\ref{eq:fractional vp}) with
 $A=a$ and with the boundary condition specified at $t=a$, should
 be reformulated as follows:
 \beq \label{eq:EL popravljene}
 \frac{\pd L}{\pd u} + {}_t^cD_b^\alpha \Big(\frac{\pd L}{\pd _aD_t^\alpha
 u}\Big) +\frac{\pd L}{\pd {}_aD_t^\alpha u}{\Big\vert}_{_{t=b}}
 \frac{1}{\Gamma(1-\alpha)}\frac{1}{(b-t)^\alpha}=0,
 \eeq
 instead of (\ref{eq:EL with Caputo}) and (\ref{eq:transversality cond}).

 We present now the Euler-Lagrange equation for (\ref{eq:fractional vp}). We
 state this as:

 \bthm
 Let $u^{\ast}\in AC([a,b])$ be an extremal of the functional
 $\cL$ in (\ref{eq:fractional vp}), whose Lagrangian $L$ satisfies
 (\ref{eq:conditions on L}). Then $u^{\ast}$
 satisfies the following Euler-Lagrange equations
 \bea
 \frac{\pd L}{\pd u}+{}_{t}^cD_{B}^{\alpha}
 \Big( \frac{\pd L}{\pd {}_{a}D_{t}^{\alpha}u}\Big)
 +\frac{\pd L}{\pd {}_aD_t^\alpha u}{\Big\vert}_{_{t=B}}
 \frac{1}{\Gamma(1-\alpha)}\frac{1}{(B-t)^\alpha}
 &=& 0,\; t\in (A,B) \label{eq:EL AB}\\
 {}_{t}D_{B}^{\alpha}\Big( \frac{\pd L}{\pd {}_{a}D_{t}^{\alpha}u}
 \Big) -{}_{t}D_{A}^{\alpha}\Big(
 \frac{\pd L}{\pd {}_{a}D_{t}^{\alpha}u}\Big)  &=& 0,\; t\in
 (a,A). \label{eq:EL aA}
 \eea
 \ethm

 \pr
 It is known that a necessary condition for a solution $u^{\ast}$
 of a fractional variational problem $\cL[u]$ defined by
 (\ref{eq:fractional vp}) is that the first variation of $\cL[u]$
 is zero at the solution $u^{\ast}$, i.e.
 \bea \label{a1}
 0 &=& \delta \mathcal{L}[u] \nonumber \\
 &=& \int_A^B \delta L(t,u(t),
 {}_{a}D_{t}^{\alpha}u)\, dt \nonumber \\
 &=& \int_A^B \frac{d}{d\eps}{\Big\vert}_{_{\eps=0}} \bigg[
 L(t,u(t)+\eps \delta u(t),{}_aD_t^\alpha (u+\eps \delta u)) -
 L(t,u(t),{}_aD_t^\alpha u)
 \bigg]\, dt \nonumber \\
 &=& \int_{A}^{B}\bigg[\frac{\partial L}{\partial u}
 \delta u(t) +\frac{\partial L}{\partial {}_{a}D_{t}^{\alpha}u}
 {}_aD_{t}^{\alpha }\delta u(t)\bigg]\, dt, \label{a1}
 \eea
 where $\delta u$ is the Lagrangian variation of $u$, i.e.
 $\delta u(a)=0$.
 Integration by parts formula (\ref{eq:frac int by parts}) gives that
 $$
 \int_{a}^{B}\frac{\partial L}{\partial {}_{a}D_{t}^{\alpha}u}
 {}_aD_{t}^{\alpha }\delta u(t) \,dt
 =
 \int_{a}^{B}\delta u(t) {}_{t}D_{B}^{\alpha }
 \Big(\frac{\partial L}{\partial {}_{a}D_{t}^{\alpha }u}\Big)\,
 dt.
 $$
 Thus we obtain
 \begin{eqnarray*}
 \int_{a}^{B}\frac{\partial L}{\partial {}_{a}D_{t}^{\alpha}u}
 {}_aD_{t}^{\alpha}\delta u(t)\, dt
 &=&
 \int_{A}^{B}\frac{\partial L}{\partial {}_{a}D_{t}^{\alpha}u}
 {}_aD_{t}^{\alpha}\delta u(t)\, dt\\
 && \qquad +\int_{a}^{A}\frac{\partial L}{\partial {}_{a}D_{t}^{\alpha}u}
 {}_aD_{t}^{\alpha}\delta u(t)\, dt \\
 &=&
 \int_{A}^{B}\delta u(t) {}_{t}D_{B}^{\alpha}\Big(
 \frac{\partial L}{\partial {}_{a}D_{t}^{\alpha}u}\Big)\, dt\\
 && \qquad +\int_{a}^{A}\delta u(t) {}_{t}D_{B}^{\alpha}
 \Big(\frac{\partial L}{\partial {}_{a}D_{t}^{\alpha}u}\Big)\,
 dt.
 \end{eqnarray*}
 From the last equality we conclude that
 \begin{eqnarray*}
 \int_{A}^{B}\frac{\partial L}{\partial {}_{a}D_{t}^{\alpha}u}
 {}_aD_{t}^{\alpha}\delta u(t)\, dt
 &=& \int_{A}^{B}\delta u(t) {}_{t}D_{B}^{\alpha}
 \Big(\frac{\partial L}{\partial _{a}D_{t}^{\alpha }u}\Big)\, dt \\
 && \qquad +\int_{a}^{A}\delta u(t) {}_{t}D_{B}^{\alpha }
 \Big(\frac{\partial L}{\partial {}_{a}D_{t}^{\alpha}u}\Big)\, dt \\
 && \qquad -\int_{a}^{A}\frac{\partial L}{\partial {}_{a}D_{t}^{\alpha}u}
 {}_aD_{t}^{\alpha }\delta u(t)\, dt\\
 &=& \int_{A}^{B}\delta u(t) {}_{t}D_{B}^{\alpha}
 \Big(\frac{\partial L}{\partial {}_{a}D_{t}^{\alpha }u}\Big)\, dt \\
 && \!\!\!+\int_{a}^{A} \bigg[ {}_{t}D_{B}^{\alpha }
 \Big(\frac{\partial L}{\partial {}_{a}D_{t}^{\alpha}u}\Big)
 -{}_{t}D_{A}^{\alpha }
 \Big(\frac{\partial L}{\partial {}_{a}D_{t}^{\alpha}u}\Big)
 \bigg]\delta u(t)\, dt.
 \end{eqnarray*}
 If we insert this into (\ref{a1}) we obtain
 $$
 0 = \int_A^B \bigg[
 \frac{\partial L}{\partial u}+{}_{t}D_{B}^{\alpha}
 \Big(\frac{\partial L}{\partial {}_{a}D_{t}^{\alpha }u}\Big)
 \bigg]\delta u(t)\, dt
 +
 \int_a^A \bigg[ {}_{t}D_{B}^{\alpha}
 \Big(\frac{\partial L}{\partial {}_{a}D_{t}^{\alpha}u}\Big)
 -{}_{t}D_{A}^{\alpha } \Big(\frac{\partial L}{\partial
 {}_{a}D_{t}^{\alpha}u}\Big)
 \bigg] \delta u(t)\, dt.
 $$
 Therefore,
 $$
 \frac{\partial L}{\partial u}+{}_{t}D_{B}^{\alpha}
 \frac{\partial L}{\partial {}_{a}D_{t}^{\alpha }u}=0,\quad
 t\in (A,B)
 $$
 and
 $$
 {}_{t}D_{B}^{\alpha}\Big(
 \frac{\partial L}{\partial {}_{a}D_{t}^{\alpha}u}\Big)
 -{}_{t}D_{A}^{\alpha }\Big(
 \frac{\partial L}{\partial {}_{a}D_{t}^{\alpha}u}\Big),
 \quad t\in (a,A)
 $$
 The claim now follows if we replace the right Riemann-Liouville
 by the right Caputo fractional derivative according to
 (\ref{eq:veza RL i Caputo-desni}) in the first equation.
 \ep

 \brem
 It is interesting to compare the Euler-Lagrange equations
 (\ref{eq:EL popravljene}) and (\ref{eq:EL AB},\ref{eq:EL aA}) in
 the case $B=b$, when $A> a$ and $A\to a$. Thus, if we let $A\to a$ in
 (\ref{eq:EL AB},\ref{eq:EL aA}) we obtain the Euler-Lagrange
 equation (\ref{eq:EL popravljene}) plus an additional condition
 $$
 {}_AI_B^{1-\alpha}\Big(\frac{\pd L}{\pd {}_aD_t^\alpha u}\Big)\equiv
 \mbox{const.}
 $$
 Indeed, since
 $$
 0={}_{t}D_{B}^{\alpha}(\pd_3 L)-{}_{t}D_{A}^{\alpha}(\pd_3 L) =
 -\frac{d}{dt}\frac{1}{\Gamma (1-\alpha)}
 \int_A^B\frac{\pd_3 L(\theta,u(\theta),
 {}_aD_\theta^\alpha u)}{(\theta-t)^\alpha}\, d\theta,
 \quad t\in(a,A),
 $$
 we obtain that ${}_AI_B^{1-\alpha} (\pd_3 L)\equiv
 \mbox{const}$.
 \erem

 \brem
 Comparing fractional calculus used in
 \cite{Jumarie06, Jumarie06a, Jumarie07} and here, one can notice the
 difference in approaches: in his work G. Jumarie used a modified
 Riemann-Liouville fractional derivative that has different properties
 (see e.g. \cite{Jumarie06a, Jumarie07}).
 \erem

 \brem
 The problem
 of formulating necessary and sufficient conditions which
 guarantee that given fractional order differential equation is
 derivable from a variational principle is still open
 (c.f. \cite{Stantilli} for inverse problems with integer
 order derivatives).
 However, it is shown in
 \cite{DreisigmeyerYoung04, DreisigmeyerYoung03} that a necessary
 condition for constructing a fractional Lagrangian from the
 given Euler-Lagrange equation is that the Euler-Lagrange equation
 involves both left and right Riemann-Liouville fractional
 derivatives. If only one of them appears then such an equation can not
 be the Euler-Lagrange equation for some Lagrangian. For example,
 according to \cite{DreisigmeyerYoung04, DreisigmeyerYoung03}
 it is not possible to construct a fractional Lagrangian for
 a linear oscillator with fractional derivative
 $u''+u+{}_aD_t^\alpha u=0$, and therefore the same holds for
 a nonlinear oscillator of the type $u''+f(u)+{}_aD_t^\alpha u=0$.
 \erem

 \section{Approximation of Euler-Lagrange equations}
 \label{sec:approximations}

 In this section we use the approximation of
 the Riemann-Liouville fractional derivative by the finite sum
 where derivatives of integer order appear, and in this way we
 analyze a fractional variational problem involving only classical
 derivatives.
 Then we examine relation between
 the Euler-Lagrange equations obtained in the process of
 approximations and the fractional Euler-Lagrange equations
 derived in the previous section.

 We will assume in the sequel that
 $L\in\cC^N([a,b]\times\R\times\R)$, at least.

 Let $(c,d)$, $-\infty<c<d<+\infty$, be an open interval in $\R$
 which contains $[a,b]$, such that for each $t\in [a,b]$ the closed
 ball $L(t,b-a)$, with center at $t$ and radius $b-a$, lies in $(c,d)$.

 For any real analytic function $f$ in $(c,d)$ we can write
 the following expansion formula:
 \beq \label{eq:sum of derivatives}
 {}_aD_t^\alpha f =\sum_{i=0}^{\infty} {\alpha\choose i}
 \frac{(t-a)^{i-\alpha}}{\Gamma (i+1-\alpha)} f^{(i)}(t),
 \quad t\in L(t,b-a)\subset (c,d),
 \eeq
 where $\ds {\alpha\choose i}=
 \frac{(-1)^{i-1}\alpha\Gamma(i-\alpha)}{\Gamma(1-\alpha)\Gamma(i+1)}$
 (cf.\ \cite[(15.4) and (1.48)]{SamkoKM}). Actually, condition
 $L(t,b-a)\subset (c,d)$ is not formulated in the literature; it
 comes from the Taylor expansion of $f(t-\tau)$ at $t$, for
 $\tau\in(a,t)$ and $t\in(a,b)$.

 Consider again the fractional variational problem
 (\ref{eq:fractional vp}). Assume that we are looking for a
 minimizer $u\in\cC^{2N}([a,b])$, for some $N\in\N$.
 We replace in the Lagrangian the left Riemann-Liouville fractional
 derivative ${}_aD_t^\alpha u$ by the finite sum of integer-valued
 derivatives as in (\ref{eq:sum of derivatives}):
 \bea
 &&\int_A^B L\Big(t, u(t), \sum_{i=0}^{N} {\alpha\choose i}
 \frac{(t-a)^{i-\alpha}}{\Gamma (i+1-\alpha)} u^{(i)}(t)\Big)\, dt
 \nonumber\\
 && \qquad\qquad =\int_A^B \bar{L} (t,u(t),u^{(1)}(t),u^{(2)}(t),
 \ldots,u^{(N)}(t))\, dt. \label{eq:N derivatives}
 \eea
 Now the Lagrangian $\bar{L}$ depends on
 $t, u$ and all (classical) derivatives of $u$ up to order $N$.
 Moreover, $\pd_3 \bar{L},\ldots,\pd_{N+2} \bar{L}\in
 \cC^{N-1}([a,b]\times\R\times\R)$, since $\pd_i\bar{L}=\pd_3 L
 {\alpha \choose i} \frac{(t-a)^{i-\alpha}}{\Gamma(i+1-\alpha)},
 i=3,...,N+2$.

 The Euler-Lagrange equation for (\ref{eq:N derivatives}) has the
 following form:
 $$
 \sum_{i=0}^{N} \Big(-\frac{d}{dt}\Big)^i \frac{\pd \bar{L}}{\pd
 u^{(i)}}=0.
 $$
 This is equivalent to
 \beq \label{eq:EL-N}
 \frac{\pd L}{\pd u}+\sum_{i=0}^{N} \Big(-\frac{d}{dt}\Big)^i
 \bigg(\pd_3 L\cdot  {\alpha\choose i}
 \frac{(t-a)^{i-\alpha}}{\Gamma (i+1-\alpha)}\bigg)=0.
 \eeq

 \brem
 The Euler-Lagrange equation (\ref{eq:EL-N}) provides
 a necessary condition when one solves the variational problem
 (\ref{eq:N derivatives}) in the class $\cC^{2N}([a,b])$,
 with the prescribed boundary condition at $A$ and $B$,
 i.e. $u(A)=A_0$ and $u(B)=B_0$, $A_0,B_0$ are fixed real numbers.
 \erem

 The question arises how (\ref{eq:EL-N}) is related to
 (\ref{eq:EL AB},\ref{eq:EL aA}). More precisely, we want to show
 that (\ref{eq:EL-N}) converges to (\ref{eq:EL AB},\ref{eq:EL aA}),
 as $N\to \infty$, in a weak sense.

 We will simplify the proof by considering the case $A=a$ and
 $B=b$. This choice reduce (\ref{eq:EL AB},\ref{eq:EL aA}) to
 (\ref{eq:EL popravljene}).

 First we prove a result which provides an expression for the
 right Riemann-Liouville fractional derivative in terms of the
 lower bound $a$, which figures in the left Riemann-Liouville fractional
 derivative.
 Such an equality holds in a weak sense, if for a test
 function space we use the space of real analytic functions as
 follows.

 Let $\mathcal{A}((c,d))$ be the space of real analytic
 functions in $(c,d)$ with the family of seminorms
 $$
 p_{[m,n]}(\vphi):= \sup_{t\in [m,n]}|\vphi(t)|, \quad\vphi\in
 \mathcal{A}((c,d)),
 $$
 where $[m,n]$ are subintervals of $(c,d)$. Every function
 $f\in\cC([a,b])$, which we extend to be zero in $(c,d)\backslash
 [a,b]$, defines an element of the dual $\mathcal{A}'((c,d))$ via
 $$
 \vphi\mapsto\langle f,\vphi \rangle =\int_a^b f(t)\vphi(t)dt,
 \quad \vphi\in\cA((c,d)).
 $$

 As usual, we say that $f$ and $g$ from
 $\mathcal{A}'((c,d))$ are equal in the weak sense if for every $\vphi \in
 \mathcal{A}((c,d))$,
 $$\langle f,\vphi\rangle = \langle g,\vphi\rangle.
 $$
 In the proposition and theorem which are to follow, we will assume
 (as in (\ref{eq:sum of derivatives})) that $L(t,b-a)\subset
 (c,d)$, for all $t\in[a,b]$.

 \bprop \label{prop:rrl-weak}
 Let $F\in\cC^\infty([a,b])$, such that $F^{(i)}(b)=0$,
 for all $i\in\N_0$, and $F\equiv 0$ in $(c,d)\backslash [a,b]$.
 Let ${}_tD_b^\alpha F$ be extended by zero in $(c,d)
 \backslash [a,b]$.
 Then:
 \begin{itemize}
 \item[(i)] For every $i\in\N$, the $(i-1)$-th derivative of
 $
 t\mapsto F(t)(t-a)^{i-\alpha}
 $
 is continuous at $t=a$ and $t=b$ and  the $i$-th derivative
 of this function, $i\in\N_0$, is integrable in $(c,d)$ and
 supported by $[a,b]$. \item[(ii)] The partial sums $S_N, N\in\N_0$,
 where
 $$
 t\mapsto S_N(t):=\left\{\ba{ll}
 \ds\sum_{i=0}^{N}
 \Big(-\frac{d}{dt}\Big)^i \bigg(F\cdot {\alpha\choose i}
 \frac{(t-a)^{i-\alpha}}{\Gamma (i+1-\alpha)}\bigg), & t\in [a,b],\\
 0, & t\in (c,d)\backslash [a,b],
 \ea \right.
 $$
 are  integrable functions in $(c,d)$ and supported by $[a,b];$
 \item[(iii)]
 \beq \label{eq:right RL weak}
 {}_tD_b^\alpha F = \sum_{i=0}^{\infty} \Big(-\frac{d}{dt}\Big)^i
 \bigg(F\cdot {\alpha\choose i}
 \frac{(t-a)^{i-\alpha}}{\Gamma (i+1-\alpha)}
 \bigg)
 \eeq
 in the weak sense.
 \end{itemize}
 \eprop

 \pr
 One can simply prove the assertions (i) and (ii) concerning the mappings
 $t\mapsto F(t)(t-a)^{i-\alpha}$ and $t\mapsto S_N(t)$, $t\in [a,b]$.

 So, let us prove the main assertion (iii).
 We have to show that
 $$
 \langle {}_tD_b^\alpha F,\vphi\rangle =
 \Big\langle \sum_{i=0}^{\infty} \Big(-\frac{d}{dt}\Big)^i
 \bigg(F\cdot {\alpha\choose i}
 \frac{(t-a)^{i-\alpha}}{\Gamma (i+1-\alpha)}
 \bigg),\vphi \Big\rangle, \quad \forall \vphi\in
 \mathcal{A}((c,d)).
 $$
 Since ${}_tD_b^\alpha F$ is continuous in $[a,b]$, it follows
 \beas
 \langle {}_tD_b^\alpha F,\vphi\rangle &=&
 \int_a^b {}_tD_b^\alpha F \vphi(t)\, dt\\
 &=& \int_a^b F(t) {}_aD_t^\alpha \vphi\, dt,
 \eeas
 where we have used fractional integration by parts
 (\ref{eq:frac int by parts}). Now, by (\ref{eq:sum of derivatives}),
 (i) and (ii), we continue
 \beas
 \langle {}_tD_b^\alpha F,\vphi\rangle
 &=& \int_a^b F(t) \sum_{i=0}^{\infty} {\alpha\choose i}
 \frac{(t-a)^{i-\alpha}}{\Gamma (i+1-\alpha)}
 \vphi^{(i)}(t)\, dt\\
 &=& \lim_{N\to \infty}\Big\langle F(t), \sum_{i=0}^{N} {\alpha\choose i}
 \frac{(t-a)^{i-\alpha}}{\Gamma (i+1-\alpha)}
 \vphi^{(i)}(t) \Big\rangle \\
 &=& \lim_{N\to \infty}\int_a^b \sum_{i=0}^{N} \Big(-\frac{d}{dt}\Big)^i
 \bigg(F(t)\cdot {\alpha\choose i}
 \frac{(t-a)^{i-\alpha}}{\Gamma (i+1-\alpha)}
 \bigg)\vphi(t)\, dt.
 \eeas
 This implies that
 $$
 \lim_{N\to\infty}
 \sum_{i=0}^{N}\int_a^b \Big(-\frac{d}{dt}\Big)^i
 \bigg(F(t)\cdot {\alpha\choose i}
 \frac{(t-a)^{i-\alpha}}{\Gamma (i+1-\alpha)}
 \bigg)
 $$
 exists in $\cA'((c,d))$ and
 \beas
 && \lim_{N\to\infty} \int_a^b \sum_{i=0}^{N} \Big(-\frac{d}{dt}\Big)^i
 \bigg(F(t)\cdot {\alpha\choose i}
 \frac{(t-a)^{i-\alpha}}{\Gamma (i+1-\alpha)}
 \bigg)\vphi(t)\, dt\\
 && \qquad\quad =\Big\langle \sum_{i=0}^{\infty} \Big(-\frac{d}{dt}\Big)^i
 \bigg(F\cdot {\alpha\choose i}
 \frac{(\cdot-a)^{i-\alpha}}{\Gamma (i+1-\alpha)}
 \bigg),\vphi \Big\rangle.
 \eeas
 This proves (\ref{eq:right RL weak}).
 \ep

 We will show in the theorem which is to follow, that the Euler-Lagrange
 equation (\ref{eq:EL-N}) converges to (\ref{eq:EL popravljene}),
 as $N\to +\infty$, in the weak sense. To shorten the notation,
 we introduce $P_N$ and $P$ for the Euler-Lagrange equations in
 (\ref{eq:EL-N}) and (\ref{eq:EL popravljene}) respectively.

 We will use the following assumptions:
 \begin{itemize}
 \item[a)] Let $u\in\cC^\infty([a,b])$ such that $u(a)=a_0$,
 $u(b)=b_0$, for fixed $a_0,b_0\in\R$, and $L_3$ (where $L_3$
 stands for $\pd_3 L$) be a function
 in $[a,b]$ defined by $t\mapsto L_3(t)= L_3(t,u(t),{}_aD_t^\alpha u)$, $t\in
 [a,b]$. Let $L_3^{(i)}(b,b_0,p)=0$, for all $i\in\N$, meaning that for $(t,s,p)
 \mapsto L_3(t,s,p)$, $t\in [a,b]$, $s,p\in\R$, the following holds:
 \begin{itemize}
 \item[(i)] $\ds \frac{\pd^i L_3}{\pd t^i}(b,b_0,p)=0, \forall p\in\R$;
 \item[(ii)] $\ds \frac{\pd^i L_3}{\pd s^i}(b,b_0,p)=0, \forall p\in\R$;
 \item[(iii)] $\ds \frac{\pd^i L_3}{\pd p^i}(b,b_0,p)=0, \forall
 p\in\R$.
 \end{itemize}
 \item[b)] Let $u\in\cC^\infty([a,b])$ such that $u^{(i)}(b)=0$,
 for all $i\in\N_0$, and $u(a)=a_0$, for fixed
 $a_0\in\R$. Let $L_3^{(i)}(b)=L_3^{(i)}(b,0,
 {}_aD_b^\alpha u)=0$, for all $i\in\N$ and for every fixed $u$,
 meaning that for $(t,s,p) \mapsto L_3(t,s,p)$, $t\in [a,b]$, $s,p\in\R$,
 the following holds:
 \begin{itemize}
 \item[(i)] $\ds \frac{\pd^i L_3}{\pd t^i}(b,0,p)=0, \forall p\in\R$;
 \item[(ii)] $\ds \frac{\pd^i L_3}{\pd p^i}(b,0,p)=0, \forall
 p\in\R$.
 \end{itemize}
 \end{itemize}
 Now we will consider the fractional variational problem
 (\ref{eq:fractional vp}) in the case a) and in the case b).

 \bthm \label{th:EL weak approx}
 Let $\cL[u]$ be a fractional variational problem
 (\ref{eq:fractional vp}) which is being
 solved in the case a) or b).
 Denote by $P$ the
 fractional Euler-Lagrange equations (\ref{eq:EL AB}),
 and by $P_N$ the Euler-Lagrange equations (\ref{eq:EL-N}), which
 correspond to the variational problem (\ref{eq:N derivatives}),
 in which the left Riemann-Liouville fractional derivative is
 approximated according to (\ref{eq:sum of derivatives}) by the
 finite sum. Then in both cases a) and b)
 $$
 P_N\to P \mbox{ in the weak sense, as } N\to \infty.
 $$
 \ethm

 \pr
 The proof of the theorem is based on Proposition \ref{prop:rrl-weak}.
 By assumptions a) and b) and the extensions of partial
 derivatives of $L$ to be zero in $(c,d)\backslash [a,b]$ we can
 apply Proposition \ref{prop:rrl-weak} with $F(t)=\pd_3 L(t,u(t),
 {}_aD_t^\alpha u)$, $t\in [a,b]$ ($F\equiv 0$ in $(c,d)\backslash [a,b]$).
 For any $\vphi\in \mathcal{A}((c,d))$ the following holds:
 \beas
 && \lim_{N\to +\infty}
 \Big\langle \frac{\pd L}{\pd u}(t,u(t),{}_aD_t^\alpha u)\\
 &&\qquad\qquad +\sum_{i=0}^{N}
 \Big(-\frac{d}{dt}\Big)^i \bigg(\pd_3 L(t,u(t),{}_aD_t^\alpha u)
 \cdot {\alpha\choose i}
 \frac{(t-a)^{i-\alpha}}{\Gamma (i+1-\alpha)}\bigg), \vphi(t)
 \Big\rangle  \\
 && \qquad\quad
 = \Big\langle \frac{\pd L}{\pd u}, \vphi \Big\rangle
 + \lim_{N\to +\infty} \Big\langle \pd_3 L,
 \sum_{i=0}^{\infty} {\alpha\choose i}
 \frac{(\cdot-a)^{i-\alpha}}{\Gamma (i+1-\alpha)}
 \vphi^{(i)} \Big\rangle\\
 && \qquad\quad
 = \Big\langle \frac{\pd L}{\pd u}, \vphi \Big\rangle
 + \langle \pd_3 L, {}_aD_t^\alpha \vphi \rangle\\
 && \qquad\quad
 = \Big\langle \frac{\pd L}{\pd u}+{}_tD_b^\alpha
 \frac{\pd L}{\pd {}_aD_t^\alpha u}, \vphi \Big\rangle.
 \eeas
 The claim now follows from (\ref{eq:veza RL i Caputo-levi}).
 \ep

 \section{Concluding remarks}

 Euler-Lagrange equations have been studied for a general
 fractional variational problem in which the lower bound in the
 variational integral differs from the lower bound in the left
 Riemann-Liouville fractional derivative which appears in the
 Lagrangian. Thus we allow for the possibility that the beginning
 of the memory of the system ($t=a$) does not coincide with the
 lower bound ($t=A$) in the Hamiltonian's action integral.
 Also, the previous results related to fractional Euler-Lagrange
 equations have been corrected and improved.

 An approximation of
 fractional derivatives in the Lagrangian has been suggested,
 resulting in a derivation of approximate Euler-Lagrange
 equations.
 Since the Leibnitz formula does not hold for ${}_aD_t^\alpha
 (f\cdot g)$, the passage from the approximated to fractional
 Euler-Lagrange equations is done by the use of weak limits
 over a specified test function space. In this way right and
 left Riemann-Liouville fractional derivatives are related to
 each other in a weak sense.

 The further research will continue towards fractional variational
 symmetries and N\"{o}ther's theorem. In this context Theorem
 \ref{th:EL weak approx} has an important role, since in a similar
 manner will be approximated the corresponding infinitesimal
 criterion as well as N\"{o}ther's theorem, which leads to a
 further analysis of variational symmetries through fractional
 calculus.



\begin{thebibliography}{10}

\bibitem{Agrawal02}
{Agrawal, O.P.}
\newblock Formulation of {E}uler-{L}agrange equations for fractional
  variational problems.
\newblock {\em J.\ Math.\ Anal.\ Appl.}, {\bf 272}:368--379, 2002.

\bibitem{Agrawal06}
{Agrawal, O.P.}
\newblock Fractional variational calculus and the transversality conditions.
\newblock {\em J.\ Phys.\ A:\ Math.\ Gen.}, {\bf 39}:10375--10384, 2006.

\bibitem{Agrawal07}
{Agrawal, O.P.}
\newblock Fractional variational calculus in terms of {R}iesz fractional
  derivatives.
\newblock {\em J.\ Phys.\ A:\ Math.\ Theor.}, {\bf 40}:6287--6303, 2007.

\bibitem{AOP07}
{Atanackovi\' c, T.M., Oparnica, Lj., Pilipovi\' c, S.}
\newblock On a nonlinear distributed order fractional differential equation.
\newblock {\em J.\ Math.\ Anal.\ Appl.}, {\bf 328}:590--608, 2007.

\bibitem{AtanackovicStankovic07}
{Atanackovi\' c, T.M., Stankovi\' c, B.}
\newblock On a class of differential equations with left and right fractional
  derivatives.
\newblock {\em Z.\ Angew.\ Math.\ Mech.}, {\bf 87}(7):537--546, 2007.

\bibitem{Bal06}
{Baleanu, D.}
\newblock Fractional {H}amiltonian analysis of irregular systems.
\newblock {\em Signal Processing}, {\bf 86}:2632--2636, 2006.

\bibitem{Bal04}
{Baleanu, D., Avkar, T.}
\newblock {L}agrangians with linear velocities within {R}iemann-{L}iouville
  fractional derivatives.
\newblock {\em Il Nuovo Cimento B}, {\bf 119}(1):73--79, 2004.

\bibitem{Bal05}
{Baleanu, D., Muslih, S.I.}
\newblock {H}amiltonian formulation of systems with linear velocities within
  {R}iemann-{L}iouville fractional derivatives.
\newblock {\em J.\ Math.\ Anal.\ Appl.}, {\bf 304}(2):599--606, 2005.

\bibitem{DreisigmeyerYoung04}
{Dreisigmeyer, D. W., Young, P. M.}
\newblock Extending {B}auer's corollary to frcational derivatives.
\newblock {\em J.\ Phys.\ A: Math. Gen.}, {\bf 37}:117--121, 2003.

\bibitem{DreisigmeyerYoung03}
{Dreisigmeyer, D. W., Young, P. M.}
\newblock Nonconservative {L}agrangian mechanics: a generalized function
  approach.
\newblock {\em J.\ Phys.\ A: Math. Gen.}, {\bf 36}:8297--8310, 2003.

\bibitem{FredericoTorres}
{Frederico, G.S.F., Torres, D.F.M.}
\newblock A formulation of {N}\"{o}ther's theorem for fractional problems of
  the calculus of variations.
\newblock {\em J.\ Math.\ Anal.\ Appl.}, {\bf 334}(2):834--846, 2007.

\bibitem{GaiesElAkrmi}
{Gaies, A., El-Akrmi, A.}
\newblock Fractional variational principle in macroscopic picture.
\newblock {\em Physica Scripta}, {\bf 70}:7--10, 2004.

\bibitem{GelfandF}
{Gelfand, I.M., Fomin, S.V.}
\newblock {\em Calculus of Variations}.
\newblock Dover publications, Inc., Mineola, New York, 2000.

\bibitem{GiaquintaH}
{Giaquinta, M., Hildebrandt, S.}
\newblock {\em Calculus of Variations I}.
\newblock Springer-Verlag, Berlin, 1996.

\bibitem{GorenfloMainardi}
{Gorenflo, R., Mainardi, F.}
\newblock Fractional calculus: Integral and differential equations of
  fractional order.
\newblock In Mainardi~F. Carpinteri, A., editor, {\em Fractals and Fractional
  Calculus in Continuum Mechanics,}, volume~{\bf 378} of {\em CISM Courses and
  Lectures}, pages 223--276. Springer-Verlag, Wien and New York, 1997.

\bibitem{Jost}
{Jost, J., Li-Jost, X.}
\newblock {\em Calculus of Variations}.
\newblock Cambridge University Press, 1998.

\bibitem{Jumarie93}
{Jumarie, G.}
\newblock Stochastic differential equations with fractional {B}rownian motion
  input.
\newblock {\em Int.\ J.\ Syst.\ Sci.}, {\bf 6}:1113--1132, 1993.

\bibitem{Jumarie95}
{Jumarie, G.}
\newblock A practical variational approach to stochastic optimal control via
  state moment equations.
\newblock {\em J.\ Franklin Inst.: Eng.\ Appl.\ Math.}, {\bf
  332(B)}(6):761--772, 1995.

\bibitem{Jumarie06}
{Jumarie, G.}
\newblock Fractionalization of the complex-valued {B}rownian motion of order
  $n$ using {R}iemann-{L}iouville derivative. {A}pplications to mathematical
  finance and stochastic mechanics.
\newblock {\em Chaos Solitons Fractals}, {\bf 28}:1285--1305, 2006.

\bibitem{Jumarie06a}
{Jumarie, G.}
\newblock Modified {R}iemann-{L}iouville derivative and fractional {T}aylor
  series of nondifferentiable functions further results.
\newblock {\em Comput.\ Math.\ Appl.}, {\bf 51}:1367--1376, 2006.

\bibitem{Jumarie07}
{Jumarie, G.}
\newblock Lagrangian mechanics of fractional order, {H}amilton-{J}acobi
  fractional pde and {T}aylor's series of nondifferential functions.
\newblock {\em Chaos Solitons Fractals}, {\bf 32}:969--987, 2007.

\bibitem{KilSriTru}
{Kilbas, A.A., Srivastava, H.M., Trujillo, J.J.}
\newblock {\em Theory and Applications of Fractional Differential Equations}.
\newblock Elsevier, Amsterdam, 2006.

\bibitem{Klimek01}
{Klimek, M.}
\newblock Fractional sequential mechanics - model with symmetric fractional
  derivative.
\newblock {\em Czech. J. Phys.}, {\bf 51}:1348--1354, 2001.

\bibitem{Lazopoulos06}
{Lazopoulos, K.A.}
\newblock Non-local continuum mechanics and fractional calculus.
\newblock {\em Mech.\ Res.\ Commun.}, {\bf 33}:753--757, 2006.

\bibitem{MillRoss}
{Miller, K.S., Ross, B.}
\newblock {\em An Introduction to the Fractional Integrals and Derivatives -
  Theory and Applications}.
\newblock John Willey \& Sons, Inc., New York, 1993.

\bibitem{Odibat07}
{Momani, S., Odibat, Z.}
\newblock Numerical comparision of methods for solving linear differential
  equations of fractional order.
\newblock {\em Chaos Solitons Fractals}, {\bf 31}(5):1248--1255, 2007.

\bibitem{Odibat06}
{Odibat, Z., Momani, S.}
\newblock Application of variational iteration method to nonlinear differential
  equations of fractional order.
\newblock {\em Int.\ J.\ Nonlinear Sci.\ Numer.\ Simul.}, {\bf 7}(1):27--34,
  2006.

\bibitem{OldSpa}
{Oldham, K.B., Spanier, J.}
\newblock {\em The Fractional Calculus}.
\newblock Academic Press, New York, 1974.

\bibitem{Baleanu&co07}
{Rabei, E.M., Nawafleh, K.I., Hijjawi, R.S., Muslih, S.I.,
Baleanu, D.}
\newblock The {H}amilton formalism with fractional derivatives.
\newblock {\em J. Math. Anal. Appl.}, {\bf 327}:891--897, 2007.

\bibitem{Riewe96}
{Riewe, F.}
\newblock Nonconservative {L}agrangian and {H}amiltonian mechanics.
\newblock {\em Phys.\ Rev. E}, {\bf 53}(2):1890--1899, 1996.

\bibitem{Riewe97}
{Riewe, F.}
\newblock Mechanics with fractional derivatives.
\newblock {\em Phys.\ Rev. E}, {\bf 55}(3):3581--3592, 1997.

\bibitem{Sagan}
{Sagan, H.}
\newblock {\em Introduction to the Calculus of Variations}.
\newblock Dover Publications, Inc., New York, 1992.

\bibitem{SamkoKM}
{Samko, S.G., Kilbas, A.A., Marichev, O.I.}
\newblock {\em Fractional Integrals and Derivatives - Theory and Applications}.
\newblock Gordon and Breach Science Publishers, Amsterdam, 1993.

\bibitem{Stantilli}
{Stantilli, R. M.}
\newblock {\em Foundations of Theoretical Mechanics I}.
\newblock Springer, Berlin, 1978.

\bibitem{vanBrunt}
{Van Brunt, B.}
\newblock {\em The Calculus of Variations}.
\newblock Springer-Verlag, New York, 2004.

\end{thebibliography}
 \end{document}